\documentclass[a4paper, BCOR=0.0mm, DIV=calc]{scrartcl}
\usepackage[T1]{fontenc}
\usepackage[utf8]{inputenc}
\usepackage[ngerman]{babel}
\usepackage[osf,sc]{mathpazo}
\usepackage{textcomp}
\usepackage{microtype}
\usepackage{amsmath, amssymb, amsfonts}
\usepackage{tabularx}
\KOMAoptions{twoside=false, twocolumn=false, headinclude=false, footinclude=false, mpinclude=false, pagesize=auto}
\recalctypearea

\setkomafont{section}{\mdseries\scshape\Large}
\setkomafont{title}{\mdseries\scshape}

\begin{document}
\title{Eine neue unendliche Reihe, die sehr stark konvergiert und die Perimetrie einer Ellipse ausdrückt\footnote{
Originaltitel: "`Nova series infinita maxime convergens perimetrum ellipsis exprimens"', erstmals publiziert in "`\textit{Novi Commentarii academiae scientiarum Petropolitanae} 18, 1774, pp. 71-84"', Nachdruck in "`\textit{Opera Omnia}: Series 1, Volume 20, pp. 357 - 370"', Eneström-Nummer E448, übersetzt von: Alexander Aycock, Textsatz: Artur Diener,  im Rahmen des Projektes "`Eulerkreis Mainz"' }}
\author{Leonhard Euler}
\date{}
\maketitle
\paragraph{§1}
Nachdem ich einst viel damit beschäftigt gewesen war, dass ich viele unendliche Reihen, mit denen die Perimetrie einer beliebigen Ellipse ausgedrückt werden würde, zu finden, hatte ich kaum vermutet, dass noch einfachere und für die Berechnung geeignetere Reihen solcher Art gefunden werden können, als die, die ich verstreut entweder in \textit{Comment. Petrop.} oder in \textit{Actis Berolin} gegeben habe.
\paragraph{§2}
Weil nun aber zufällig meine Überlegungen in dieselbe Sache hineingestolpert sind, hat sich mir eine andere, und, wenn ich mich nicht täusche, um Vieles einfachere und gefälligere Reihe offenbart, deren Untersuchungen ich im Geiste so angestellt habe: Ich betrachte natürlich den elliptischen Quadranten $ACB$ (Fig. 1, p.$359$), dessen Halbachsen $CA=a$, $CB=b$ seien; die diesen parallelen Koordinaten nenne man $CP=x$ und $PM=y$, sodass man aus der Gestalt der Ellipse diese Gleichung hat:
\[
	bbx^2 + aay^2 = aa\cdot bb
\]
oder
\[
	\frac{x^2}{a^2} + \frac{y^2}{b^2} = 1,
\]
aus welcher ich auf einzigartige Weise die Länge des ganzen Bogen $AMB$ oder den $4$.Teil der Perimetrie bestimme.
\paragraph{§3}
Weil also
\[
	\frac{x^2}{a^2} + \frac{y^2}{b^2} = 1
\]
sein muss, führe ich eine neue Variable $z$ in die Rechnung ein, indem ich
\[
	\frac{x^2}{a^2} = \frac{1+z}{2}
\]
setze, dass
\[
	\frac{y^2}{b^2} = \frac{1-z}{2}
\]
wird, woher
\[
	x = a\sqrt{\frac{1+z}{2}} \quad \text{und} \quad y = b\sqrt{\frac{1-z}{2}}
\]
hervorgeht und daher, indem man differenziert,
\[
	\mathrm{d}x = \frac{a\mathrm{d}z}{2\sqrt{2(1+z)}} \quad \text{und} \quad \mathrm{d}y = \frac{-b\mathrm{d}z}{2\sqrt{2(1-z)}};
\]
daraus berechnen wir, wenn wir den Bogen $BM=s$ nennen, sofort
\[
	\mathrm{d}s^2 = \mathrm{d}x^2 + \mathrm{d}y^2 = \frac{a^2 \mathrm{d}z^2}{8(1+z)} + \frac{b^2\mathrm{d}z^2}{8(1-z)}
\]
oder
\[
	\mathrm{d}s^2 = \frac{\mathrm{d}z^2}{8}\left(\frac{a^2}{1+z} + \frac{b^2}{1-z}\right) = \frac{\mathrm{d}z^2 (a^2 + b^2 - (a^2-b^2)z)}{8(1-z^2)}
\]
und daher durch Integrieren
\[
	s = \frac{1}{2\sqrt{2}}\int\mathrm{d}z\sqrt{\frac{a^2 + b^2-(a^2-b^2)z}{1-z^2}},
\]
nachdem das Integral so genommen wurde, dass es für $x=0$ oder $z=-1$ gesetzt verschwindet; dann aber erstrecke man das Integral bis hin zur Grenze $x=a$, wo $z=+1$ wird, und so wird man gesuchten Quadranten $AMB$ erhalten.
\paragraph{§4}
Damit wir diese Formel als der Behandlung zugänglich ausgeben, wollen wir der Kürze wegen
\[
	a^2 + b^2 = c^2 \quad \text{und} \quad \frac{a^2-b^2}{a^2+b^2} = n
\]
setzen. Auf diese Weise folgern wir nämlich
\[
	s = \frac{c}{2\sqrt{2}}\int\mathrm{d}z\frac{\sqrt{1-nz}}{\sqrt{1-z^2}},
\]
wo wir die obere Wurzel auf gewohnte Weise in eine Reihe verwandeln wollen:
\[
\sqrt{1-nz} = 1 - \tfrac{1}{2}nz - \tfrac{1\cdot 1}{2\cdot 4}n^2z^2 - \tfrac{1\cdot 1\cdot 3}{2\cdot 4\cdot 6}n^3z^3 - \tfrac{1\cdot 1\cdot 3\cdot 5}{2\cdot 4\cdot 6\cdot 8}n^4z^4 - \tfrac{1\cdot 1\cdot 3\cdot 5\cdot 7}{2\cdot 4\cdot 6\cdot 8\cdot 10}n^5z^5 - \mathrm{etc},
\]
welche einzelnen Terme uns zur einzigartigen Integration führen; und die beiden ersten werden freilich nach gegebenen Gesetz integriert, dass sie natürlich für $z=-1$ genommen verschwinden,
\[
	\int\frac{\mathrm{d}z}{\sqrt{1-z^2}} = \arcsin{z} - \arcsin{(-1)} = \arcsin{z}+\tfrac{1}{2}\pi
\]
und
\[
	\int\frac{z\mathrm{d}z}{\sqrt{1-z^2}} = -\sqrt{1-z^2} + 0
\]
geben; daher wird also, wenn wir $z=+1$ nehmen
\[
	\int\frac{\mathrm{d}z}{\sqrt{1-z^2}} = \pi \quad \text{und} \quad \int\frac{z\mathrm{d}z}{\sqrt{1-z^2}} = 0
\]
hervorgehen.
\paragraph{§5}
Für die übrigen Terme wollen wir die gewöhnliche allgemeine Reduktion betrachten:
\[
	\int\frac{z^{\lambda + 2}\mathrm{d}z}{\sqrt{1-z^2}} = A\int\frac{z^{\lambda}\mathrm{d}z}{\sqrt{1-z^2}} + Bz^{\lambda +1}\sqrt{1-z^2},
\]
wo
\[
	A = \frac{\lambda + 1}{\lambda + 2} \quad \text{und} \quad B = \frac{-1}{\lambda +2}
\]
sein muss, sodass
\[
	\int\frac{z^{\lambda +2}\mathrm{d}z}{\sqrt{1-z^2}} = \frac{\lambda +1}{\lambda +2}\int\frac{z^{\lambda}\mathrm{d}z}{\sqrt{1-z^2}} - \frac{1}{\lambda +2}z^{\lambda +1}\sqrt{1-z^2}
\]
ist, wo wir die Konstante nicht hinzugefügt haben, weil diese Formel schon für $z=-1$ genommen verschwindet; daher wird man, wenn man gleich $z=+1$ setzt,
\[
	\int\frac{z^{\lambda +2}\mathrm{d}z}{\sqrt{1-z^2}} = \frac{\lambda +1}{\lambda +2}\int\frac{z^{\lambda}\mathrm{d}z}{\sqrt{1-z^2}}
\]
erhalten.
\paragraph{§6}
Aus dieser Reduktion ist sofort klar, dass alle aus ungeraden Potenzen von $z$ zu entstehenden Integrale per se verschwinden; für die geraden Potenzen aber erhalten wir für unseren Zweck
\begin{align*}
\int\frac{\mathrm{d}z}{\sqrt{1-z^2}} = \pi ,& \qquad \int\frac{z^2\mathrm{d}z}{\sqrt{1-z^2}} = \tfrac{1}{2}\pi \\
\int\frac{z^4\mathrm{d}z}{\sqrt{1-z^2}} = \tfrac{1\cdot 3}{2\cdot 4}\pi ,& \qquad \int\frac{z^6\mathrm{d}z}{\sqrt{1-z^2}} = \tfrac{1\cdot 3\cdot 5}{2\cdot 4\cdot 6}\pi \\
\mathrm{etc.}
\end{align*}
\paragraph{§7}
Nachdem also diese Werte eingesetzt wurden, wird die Länge des Ellipsenquadranten berechnet, dass
\[
	AMB = \frac{c\pi}{2\sqrt{2}} \begin{cases}
	1 - \frac{1\cdot 1}{2\cdot 4}n^2 - \frac{1\cdot 1\cdot 3\cdot 5}{2\cdot 4\cdot 6\cdot 8}n^4\cdot\frac{1\cdot 3}{2\cdot 4} \\
	-\frac{1\cdot 1\cdot 3\cdot 5\cdot 7\cdot 9}{2\cdot 4\cdot 6\cdot 8\cdot 10\cdot 12}n^6\cdot\frac{1\cdot 3\cdot 5}{2\cdot 4\cdot 6} - \mathrm{etc}
	\end{cases}
\]
sein wird. Für diese Form aber wollen wir solange der Kürze wegen
\[
	AMB = \frac{c\pi}{2\sqrt{2}} \lbrace 1-\alpha n^2 - \beta n^4 - \gamma n^6 - \delta n^8 - \varepsilon n^{10} - \mathrm{etc} \rbrace
\]
schreiben, welche Koeffizienten auf folgende Weise kürzer ausgedrückt werden können:
\[
	\alpha = \frac{1\cdot 1}{2\cdot 4}\cdot\frac{1}{2} = \frac{1\cdot 1}{4\cdot 4}, \quad \frac{\beta}{\alpha} = \frac{3\cdot 5}{8\cdot 8},\quad \frac{\gamma}{\beta} = \frac{7\cdot 9}{12\cdot 12},\quad \frac{\delta}{\gamma} = \frac{11\cdot 13}{16\cdot 16}, \quad \mathrm{etc}
\]
\paragraph{§8}
Weil also die gefundenen Koeffizienten eine so einfache wie außergewöhnliche Reihe festsetzen, scheint dieser Ausdruck, den wir gefunden haben, jedenfalls besonders der Aufmerksamkeit würdig, weil die Terme sehr stark konvergieren und das natürlch für alle Ellipsen, deshalb weil $\frac{a^2-b^2}{a^2+b^2} = n$ immer ein Bruch kleiner als die Einheit ist. Wir werden natürlich
\[
	AMB = \frac{c\pi}{2\sqrt{2}}\begin{cases}
	1 - \frac{1\cdot 1}{4\cdot 4}n^2 - \frac{1\cdot 1}{4\cdot 4}\cdot\frac{3\cdot 5}{8\cdot 8}n^4 - \frac{1\cdot 1}{4\cdot 4}\cdot\frac{3\cdot 5}{8\cdot 8}\cdot\frac{7\cdot 9}{12\cdot 12}n^6 \\
	- \frac{1\cdot 1}{4\cdot 4}\cdot\frac{3\cdot 5}{8\cdot 8}\cdot \frac{7\cdot 9}{12\cdot 12}\cdot\frac{11\cdot 13}{16\cdot 16}n^8 - \mathrm{etc}
	\end{cases}
\]
haben.
\paragraph{§9}
Wir betrachten daher den Fall, in welchem unsere Ellipse ein Kreis des Radius gleich $a$ wird; dann nämlich wird $b=a$ sein, daher $c = a\sqrt{2}$ und $n=0$, woraus der Kreisquadrant hervorgeht, wie freilich all bekannt ist, gleich $\frac{1}{2}\pi a$.
\paragraph{§10}
Darauf aber taucht auch ein höchst bemerkenswerter Fall auf, in dem die Halbachse $CB = b = 0$ ist; dann nämlich wird der Ellipsenquadrant $AMB$ der Halbachse $CA = a$ selbst gleich, aber für unsere Formel wird $c=a$ und $n=1$ sein, nach Einsetzen welcher Werte wir die folgende Gleichung erhalten
\[
	a = \frac{\pi a}{2\sqrt{2}}\left( 1 - \frac{1\cdot 1}{4\cdot 4} - \frac{1\cdot 1}{4\cdot 4}\cdot \frac{3\cdot 5}{8\cdot 8} - \frac{1\cdot 1}{4\cdot 4}\cdot \frac{3\cdot 5}{8\cdot 8}\cdot \frac{7\cdot 9}{12\cdot 12} - \mathrm{etc} \right),
\]
welches genauer jener Fall selbst ist, in dem unsere Reihe möglich wenigst konvergent ist, und der deshalb unsere Aufmerksamkeit umso mehr verdient, weil die Summe dieser genau angegebenen werden kann, weil
\[
	1 - \frac{1\cdot 1}{4\cdot 4} - \frac{1\cdot 1}{4\cdot 4}\cdot \frac{3\cdot 5}{8\cdot 8}~~\mathrm{etc} \quad \text{ins Unendliche} = \frac{2\sqrt{2}}{\pi}
\]
ist.
\paragraph{§10[a]}
Wenn es irgendjemanden beliebt über diese Reise numerische Rechnungen anzustellen, so wollen wir hier die Werte der Buchstaben $\alpha$, $\beta$, $\gamma$, etc in Dezimalbrüchen beifügen, die sich so verhalten: 
\begin{align*}
	\alpha &= 0,0625000 \\
	\beta  &= 0,0146484 \\
	\gamma &= 0,0064087 \\
	\delta &= 0,0035798 \\
	\varepsilon &= 0,0022821 \\
	\zeta &= 0,0015808 \\
	\mathrm{etc},
\end{align*}
nachdem aber nur soweit fortgesetzt wurde, geht
\[
	1 - \alpha - \beta - \gamma - \delta - \varepsilon - \zeta = 0,9090002
\]
hervor; in der Tat aber findet man $\frac{2\sqrt{2}}{\pi} = 0,9003200$; daher sehen wir, dass die Summe aller folgenden Buchstaben $\eta$, $\vartheta$, $i$, $x$, etc $0,0086802$ ergeben muss.
\paragraph{§11}
Im Übrigen wird es für die numerische Berechnung nicht wenig förderlich sein bemerkt zu haben, dass unsere Koeffizienten auch auf die folgende Weise gefälliger ausgedrückt werden können:
\begin{align*}
	\alpha &= \frac{1}{16} \\
	\beta &= \frac{1}{64}\cdot\frac{15}{16} \\
	\gamma &= \frac{1}{144}\cdot\frac{15}{16}\cdot\frac{63}{64} \\
	\delta &= \frac{1}{256}\cdot\frac{15}{16}\cdot\frac{63}{64}\cdot\frac{143}{144} \\
	\varepsilon &= \frac{1}{400}\cdot\frac{15}{16}\cdot\frac{63}{64}\cdot\frac{143}{144}\cdot\frac{255}{256} \\
	\mathrm{etc.}
\end{align*}
\paragraph{§12}
Bei Gelegenheit dieser Reihe, die wir gefunden haben, wird es der Mühe Wert sein, deren Summe von der letzten aus zu untersuchen, das kann auf zweifache Weise gemacht werden; die erste Art, die ich schon einst vorgelegt habe und darauf sehr oft zum Nutzen angewendet habe, führt uns zu einer Differentialgleichung, deren Integral durch die vorgelegte Reihe selbst ausgedrückt werde. Damit nun diese Methode leichter gehandhabt werden kann, wollen wir $n=2v$ setzen, sodass die zu summierende Reihe
\[
	s = 1 - \frac{1\cdot 1}{2\cdot 2}v^2 - \frac{1\cdot 1}{2\cdot 2}\cdot\frac{3\cdot 5}{4\cdot 4}v^4 - \frac{1\cdot 1}{2\cdot 2}\cdot\frac{3\cdot 5}{4\cdot 4}\cdot\frac{7\cdot 9}{6\cdot 6}v^6 - \mathrm{etc}
\]
wird.
\paragraph{§13}
Wir wollen diese Reihe differenzieren, dann aber wiederum mit $v$ multiplizieren, dass
\[
	\frac{v\mathrm{d}s}{\mathrm{d}v} = -\frac{1\cdot 1}{2}v^2 - \frac{1\cdot 1}{2\cdot 2}\cdot\frac{3\cdot 5}{4}v^4 - \frac{1\cdot 1}{2\cdot 2}\cdot\frac{3\cdot 5}{4\cdot 4}\cdot\frac{7\cdot 9}{6}v^6 ~~\mathrm{etc}
\]
hervorgeht, welche erneut differenziert
\[
	\frac{\mathrm{d}\cdot v\mathrm{d}s}{\mathrm{d}v^2} = -1\cdot 1\cdot v - \frac{1\cdot 1}{2\cdot 2}\cdot 3\cdot 5 \cdot v^3 - \frac{1\cdot 1}{2\cdot 2}\cdot\frac{3\cdot 5}{4\cdot 4}\cdot 7\cdot 9\cdot v^5 ~~ \mathrm{etc}
\]
liefert; auf diese Weise haben wir natürlich aus den einzelnen Nennern $2$ Faktoren entfernt.
\paragraph{§14}
Nun aber wollen wir erneut mithilfe der Differentation die Nenner um $2$ neue Faktoren vermehren; für dieses Ziel wollen wir die erste mit $\sqrt{v}$ multiplizierte Gleichung differenzieren und es wird
\[
	\frac{2\mathrm{d}\cdot s\sqrt{v}}{\mathrm{d}v} = +v^{-\frac{1}{2}} - \frac{1\cdot 1}{2\cdot 2}5v^{\frac{3}{2}} - \frac{1\cdot 1}{2\cdot 2}\cdot\frac{3\cdot 5}{4\cdot 4}9v^{\frac{7}{2}} - \frac{1\cdot 1}{2\cdot 2}\cdot\frac{3\cdot 5}{4\cdot 4}\cdot\frac{7\cdot 9}{6\cdot 6}13v^{\frac{11}{2}} ~~ \mathrm{etc}
\]
hervorgehen; man differenziere diese erneut und es wird, indem man wiederum mit $2$ multipliziert,
\[
	\frac{4\mathrm{d}\mathrm{d}s\sqrt{v}}{\mathrm{d}v^2} = -v^{-\frac{3}{2}} - \frac{1\cdot 1}{2\cdot 2}\cdot 3\cdot 5\cdot v^{\frac{1}{2}} - \frac{1\cdot 1}{2\cdot 2}\cdot\frac{3\cdot 5}{4\cdot 4}\cdot 7\cdot 9\cdot v^{\frac{5}{2}} ~~ \mathrm{etc}
\]
sein, die mit $v^{\frac{5}{2}}$ multipliziert
\[
	\frac{4v^{\frac{5}{2}}\mathrm{d}\mathrm{d}s\sqrt{v}}{\mathrm{d}v^2} = -v - \frac{1\cdot 1}{2\cdot 2}\cdot 3\cdot 5\cdot v^3 - \frac{1\cdot 1}{2\cdot 2}\cdot\frac{3\cdot 5}{4\cdot 4}\cdot 7\cdot 9\cdot v^5 ~~ \mathrm{etc}
\]
erzeugt; oben aber haben wir schon
\[
	\frac{\mathrm{d}\cdot v\mathrm{d}s}{\mathrm{d} v^2} = -v - \frac{1\cdot 1}{2\cdot 2}\cdot 3\cdot 5 \cdot v^3 - \frac{1\cdot 1}{2\cdot 2}\cdot\frac{3\cdot 5}{4\cdot 4}\cdot 7\cdot 9 \cdot v^5 ~~ \mathrm{etc}
\]
gefunden; weil die Reihen gleich sind, folgern wir daher diese Gleichung
\[
	4v^{\frac{5}{2}}\mathrm{d}\mathrm{d}\cdot s\sqrt{v} = \mathrm{d}v \cdot \mathrm{d}s
\]
welche Gleichung die Relation der gesuchten Summe zur Variablen $v$ enthält.
\paragraph{§15}
Diese Gleichung also wird entwickelt ein Differential zweiter Ordnung, nachdem nämlich das Zeichen $\mathrm{d}v$ als Konstante genommen wurde, wird wegen
\begin{align*}
	\mathrm{d}\cdot s\sqrt{v} &= \mathrm{d}s\sqrt{v} + \frac{s\mathrm{d}v}{2\sqrt{v}} \\
	\mathrm{d}\mathrm{d}s\sqrt{v} &= \mathrm{d}\mathrm{d}s\sqrt{v} + \frac{\mathrm{d}v\mathrm{d}s}{\sqrt{v}} - \frac{s\mathrm{d}v^2}{4v\sqrt{v}}
\end{align*}
sein, also
\[
	4v^{\frac{5}{2}}\mathrm{d}\mathrm{d}s\sqrt{v} = 4v^3\mathrm{d}\mathrm{d}s + 4v^2\mathrm{d}v\mathrm{d}s - sv\mathrm{d}v^2;
\]
dann aber wird man wegen $\mathrm{d}\cdot v\mathrm{d}s = v\mathrm{d}\mathrm{d}s + \mathrm{d}v\mathrm{d}s$ diese Gleichung haben
\[
	v\mathrm{d}\mathrm{d}s(1-4v^2) + \mathrm{d}v\mathrm{d}s(1-4v^2) + sv\mathrm{d}v^2 = 0
\]
oder
\[
	v\mathrm{d}\mathrm{d}s + \mathrm{d}v\mathrm{d}s + \frac{sv\mathrm{d}v^2}{1-4v^2} = 0
\]
\paragraph{§16}
Die Konstruktion dieser Differentialgleichung zweiter Ordnung also ist in unserer Macht; es sei nämlich sie die Ellipse, deren Halbachsen $a$ und $b$ sein mögen und der vierte Teil ihrer Peripherie gleich $q=AMB$; dann aber nehme man
\[
	c = \sqrt{a^2 + b^2} \quad \text{und} \quad \frac{a^2 - b^2}{a^2 + b^2} = n = 2v,
\]
daher, weil
\[
	q = \frac{\pi c}{2\sqrt{2}}s
\]
ist, wird
\[
	s = \frac{2q\sqrt{2}}{\pi c}
\]
werden. Schon wird wegen $a^2 + b^2 = c^2$ und $a^2 - b^2 = 2c^2v$
\[
	a^2 = \frac{c^2(1+2v)}{2} \quad \text{und} \quad b^2 = \frac{c^2(1-2v)}{2}
\]
sein. Deshalb wird unsere Konstruktion so beschaffen sein: Nachdem als Halbachsen der Ellipse
\[
	a = c\sqrt{\frac{1+2v}{2}} \quad \text{und} \quad b = c\sqrt{\frac{1-2v}{2}}
\]
genommen wurden, sei $q$ der vierte Teil der Perimetrie dieser Ellipse und es wird für die Auflösung unserer Gleichung
\[
	s = \frac{2q\sqrt{2}}{\pi c}
\]
sein. Diese Gleichung, wenn wir $s = \frac{z}{\sqrt{v}}$ setzen, wird diese einfachere Form annehmen
\[
	\mathrm{d}\mathrm{d}z + \frac{z\cdot \mathrm{d} v^2}{4v^2(1-4v^2)} = 0,
\]
für welche $z = \frac{2q\sqrt{2v}}{\pi c}$ sein wird.
\paragraph{§17}
Diese Gleichung wird weiter zu einen Differential erster Ordnung zurückgeführt werden, indem man $z = e^{\int t\mathrm{d}v}$ setzt; dann wird nämlich
\[
	\mathrm{d}t + t^2\mathrm{d}v + \frac{\mathrm{d}v}{4v^2(1-4v^2)} = 0
\]
resultieren, wenn sich daher $t$ durch $v$ bestimmen ließe, sodass das Integral $\int t \mathrm{d}v$ bekannt werden würde, würde $z = e^{\int t\mathrm{d}v}$ sein.
\paragraph{§18}
Dies war die erste Art aus der vorgelegten unendlichen Reihe ihre Summe zu finden, wo wir natürlich anstelle der konstanten Zahl die variable Größe $v$ eingeführt haben; auf eine andere Art daselbe zu leisten, von welche schon verstreut viele Beispiele auftauchen, wird eine solche konstante Größe $n$ zurückgelassen; wir setzen aber $n=2m$, sodass unsere zu summierende Reihe
\[
	1 - \frac{1\cdot 1}{2\cdot 2}m^2 - \frac{1\cdot 1}{2\cdot 2}\cdot\frac{3\cdot 5}{4\cdot 4}m^4 - \frac{1\cdot 1}{2\cdot 2}\cdot\frac{3\cdot 5}{4\cdot 4}\cdot\frac{7\cdot 9}{6\cdot 6}m^6 - \mathrm{etc}
\]
ist.
\paragraph{§19}
Nun wollen wir aber ansetzen, dass
\[
	s = \int\mathrm{d} z \sqrt[4]{1-2m^2p}
\]
ist, nachdem natürlich nach Ausführung der Integration der variablen Größe $z$ ein bestimmter gewisser Wert zugeteilt worden war; wir wollen aber den Buchstaben $p$ auch als variabel betrachten; von welcher Art diese Funktion von $z$ genommen werden muss, dass diese Integration zu unserer unendlichen Reihe führt, werden wir auf die folgendene Weise untersuchen.
\paragraph{§20}
Nachdem aber die irrationale Formel $(1-2m^2p)^{\frac{1}{4}}$ in diese unendliche Reihe
\[	
	1 - \frac{1}{2}m^2p - \frac{1\cdot 3}{2\cdot 4}m^4p^2 - \frac{1\cdot 3\cdot 7}{2\cdot 4\cdot 6}m^6p^3 - \mathrm{etc}
\]
entwickelt wurde, wird die Größe $s$ durch die folgende Reihe der Integralformeln betimmt werden
\[
	s = z - \frac{1}{2}m^2\int p\mathrm{d}z - \frac{1\cdot 3}{2\cdot 4}m^4\int p^2\mathrm{d}z - \frac{1\cdot 3\cdot 7}{2\cdot 4\cdot 6}m^6\int p^3\mathrm{d}z ~~ \mathrm{etc}.
\]
Nun wollen wir aber festsetzen, wenn nach den einzelnen Integrationen der Variablen $z$ ein bestimmter Wert zugeteilt wird, dass dann
\begin{align*}
\int p\mathrm{d}z = \frac{1}{2}z, & \qquad \int p^2\mathrm{d}z = \frac{5}{4}\int p \mathrm{d}z\\
\int p^3\mathrm{d}z = \frac{9}{6}\int p^2 \mathrm{d}z, & \qquad \int p^4\mathrm{d}z = \frac{13}{8}\int p^3\mathrm{d}z\\
\mathrm{etc}
\end{align*}
sein wird; denn so wird nämlich
\[
	s = z\left( 1-\frac{1\cdot 1}{2\cdot 2}m^2 - \frac{1\cdot 1}{2\cdot 2}\cdot\frac{3\cdot 5}{4\cdot 4}m^4 - \frac{1\cdot 1}{2\cdot 2}\cdot\frac{3\cdot 5}{4\cdot 4}\cdot\frac{7\cdot 9}{6\cdot 6}m^6 ~~ \mathrm{etc}\right)
\]
werden, welche unsere vorgelegte Reihe selbst ist.
\paragraph{§21}
Nun also geht die ganze Frage darauf zurück, eine Funktion welcher Art von $z$ für $p$ genommen werden muss, dass jenes festgesetzte Verhältnis der Integrale, während natürlich der Variable $z$ ein bstimmter Wert zugeteilt wird, erhalten wird; diese Relation aber wird im Allgemeinen so ausgedrückt:
\[
	\int p^{\lambda}\mathrm{d}z = \frac{4\lambda -3}{2\lambda}\int p^{\lambda -1}\mathrm{d}z;
\]
wir wollen also, nachdem die Integrale noch unbestimmt genommen wurden, festsetzen, dass
\[
	\int p^{\lambda}\mathrm{d}z = \frac{4\lambda -3}{2\lambda}\int p^{\lambda -1}\mathrm{d}z + \frac{p^{\lambda}Q}{2\lambda}
\]
sein wird; nachdem also differenziert worden ist, wird
\[
	p^{\lambda}\mathrm{d}z = \frac{4\lambda -3}{2\lambda}p^{\lambda -1}\mathrm{d}z + \frac{1}{2}p^{\lambda-1}Q\mathrm{d}p + \frac{p^{\lambda}\mathrm{d}Q}{2\lambda}
\]
hervorgehen, welche durch $p^{\lambda -1}$ geteilt und mit $2\lambda$ multipliziert
\[
	2\lambda p\mathrm{d}z = (4\lambda -3)\mathrm{d}z + \lambda Q\mathrm{d}p + p\mathrm{d}Q
\]
liefert, und weil diese Gleichung bestehen muss, was auch immer $\lambda$ ist, liefert das uns diese zwei Gleichungen
\[
	2p\mathrm{d}z - 4\mathrm{d}z - Q\mathrm{d}p = 0,\quad -3\mathrm{d}z + p\mathrm{d}Q = 0,
\] 
aus denen sich jeder der beiden Funktionen $p$ und $Q$ bestimmen lassen wird.
\paragraph{§22}
Es ist hier aber völlig egal, ob $p$ und $Q$ Funktionen von $z$ sind oder $z$ und $Q$ von $p$, solange deren Beziehung untereinander festgesetzt wird; aus der ersten aber haben wir sofort:
\[
	\mathrm{d}z = \frac{1}{3}p\mathrm{d}Q,
\]
welcher Wert in die erste eingesetzt
\[
	\frac{2}{3}(p-2)p\mathrm{d}Q - Q\mathrm{d}p = 0
\]
liefert, aus welcher
\[
	\frac{dQ}{Q} = \frac{3\mathrm{d}p}{2p(p-2)} = -\frac{3\mathrm{d}p}{4p} + \frac{3}{4}\frac{\mathrm{d}p}{(p-2)}
\]
wird, woher durch Integrieren
\[
	\log{Q} = -\frac{3}{4}\log{P} + \frac{3}{4}\log{(p-2)} = +\frac{3}{4}\log{\frac{p-2}{p}}
\]
entsteht, woher
\[
	Q = 2\left(\frac{p-2}{p}\right)^{\frac{3}{4}}
\]
wird; dann aber, weil aus der ersten Gleichung
\[
	\mathrm{d}z = \frac{Q\mathrm{d}p}{2(p-2)}
\]
ist, wird daher
\[
	\mathrm{d}z = \frac{\mathrm{d}p}{p^{\frac{3}{4}}(p-2)^{\frac{1}{4}}} = \frac{\mathrm{d}p}{\sqrt[4]{p^3(p-2)}}
\]
sein. Nun aber muss besonders bemerkt werden, dass für jede der beiden Integrationsgrenzen die dort hinzugefügte algebraische Gleichung
\[
	p^{\lambda}Q = 2p^{\lambda - \frac{3}{4}}(p-2)^{\frac{3}{4}}
\]
verschwindet, und so ist klar, dass die Integrationsgrenzen $p=0$ und $p=2$ gesetzt werden müssen.
\paragraph{§23}
Siehe also unsere eingangs eingeführte Integralformel; auf diese Weise dargestellt
\[
	s = \int\frac{\mathrm{d}p\sqrt[4]{1-2m^2p}}{\sqrt[4]{p^3(p-2)}},
\]
weil daher
\[
	z = \int\frac{\mathrm{d}p}{\sqrt[4]{p^3(p-2)}}
\]
ist, wird unsere vorgelegte Reihe selbst
\[
	1 - \frac{1\cdot 1}{2\cdot 2}m^2 - \frac{1\cdot 1}{2\cdot 2}\cdot\frac{3\cdot 5}{4\cdot 4}m^4 ~~ \mathrm{etc}
\]
dem Bruch $\frac{s}{z}$ gleichwerden, nachdem natürlich diese Integrale so genommen worden waren, dass sie für $p=0$ gesetzt verschwinden, dann aber setze man $p=2$; deshalb sollten jene zwei Integralformeln so ausgedrückt werden:
\[
	s =\int\frac{\mathrm{d}p\sqrt[4]{1-2m^2p}}{\sqrt[4]{p^3(2-p)}} \quad \text{und} \quad z = \int\frac{\mathrm{d}p}{\sqrt[4]{p^3(2-p)}}
\]
\paragraph{§24}
Aus diesen also kann unsere oben gefundene Reihe
\[
	1 - \frac{1\cdot 1}{4\cdot 4}n^2 - \frac{1\cdot 1}{4\cdot 4}\cdot\frac{3\cdot 5}{8\cdot 8}n^4 ~~ \mathrm{etc},
\]
deren Summe wir schon gesehen haben $\frac{2q\sqrt{2}}{\pi c}$ zu sein, auch auf diese Weise durch zwei Integralformeln dargestellt werden, die, nachdem die leichte Änderung $p=2q$ gemacht wurde, sein werden wie folgt: die, die den Zähler festsetzt, wird
\[
	s = \int\frac{\mathrm{d}r\sqrt[4]{1-nnr}}{\sqrt[4]{r^3(1-r)}}
\]
sein, die andere aber, die den Nenner festsetzt,
\[
	z = \int\frac{\mathrm{d}r}{\sqrt[4]{r^3(1-r)}};
\]
dieser Bruch selbst aber wird unsere Reihe beschaffen; nun aber sind die Integrationsgrenzen $r=0$ und $r=1$.
\paragraph{§25}
Noch kürzer können diese Formeln transformiert werden, indem man $r=t^4$ nimmt; dann nämlich werden die beiden Integralformeln
\[
	s = \int\frac{\mathrm{d}t\sqrt[4]{1-n^2t^4}}{\sqrt[4]{1-t^4}} \quad \text{und} \quad z = \int\frac{\mathrm{d}t}{\sqrt[4]{1-t^4}},
\]
während die Integrationsgrenzen noch immer $t=0$ und $t=1$ werden; nachdem das bemerkt worden ist, wird der Bruch $\frac{s}{z}$ unserer Reihe gleich werden oder es wird
\[
	\frac{s}{z} = \frac{2q\sqrt{2}}{\pi c}
\]
sein, wo $q$ den vierten Teil der Peripherie einer Ellipse bezeichnet, deren Halbachsen
\[
	c\sqrt{\frac{1+n}{2}} \quad \text{und} \quad c\sqrt{\frac{1-n}{2}}
\]
sind.
\paragraph{§26}
Daher wird im Fall $n=0$ natürlich $\frac{s}{z} = 1$, im Fall $n=1$ aber wegen $s=t=1$ wird
\[
	\frac{1}{z} = \frac{2\sqrt{2}}{\pi} \quad \text{oder} \quad z = \int\frac{\mathrm{d}t}{\sqrt[4]{1-t^4}} = \frac{\pi}{2\sqrt{2}}
\]
werden, was freilich schon anderswoher bekannt ist.
\end{document}